\documentclass[a4paper,11pt]{amsart}

\usepackage{amsmath}

\def\ei{\frac{e_i}{r}}

\def\pa{\partial}
\def\i{\lrcorner}

\def\cd{\cdot}

\def\a{\alpha}

\def\Ric{{\rm Ric}}

\def\la{\langle}
\def\ra{\rangle}
\def\.{\cdot}
\def\O{\Omega}
\def\dr{{\partial_r}}
\def\n{\nabla}

\def\s{\sigma}

\def\beq{\begin{eqnarray*}}
\def\eeq{\end{eqnarray*}}

\def\x{\times}
\def\f{\varphi}

\def\o{\omega}

\def\d{\delta}

\def\nn{\bar\nabla}
\def\gg{\bar g}
\def\L{\Lambda}

\def\r{\end{proof}}

\def \RM{\mathbb{R}}

\newtheorem{ede}{Definition}[section]
\newtheorem{epr}[ede]{Proposition}
\newtheorem{ath}[ede]{Theorem}
\newtheorem{elem}[ede]{Lemma}

\newtheorem{exa}[ede]{Example}

\title{The odd--dimensional Goldberg Conjecture}

\author{Vestislav Apostolov, Tedi Dr\u aghici and Andrei Moroianu}

\thanks{The first author was supported in part by NSERC grant OGP0023879}

\address{Vestislav Apostolov\\ D{\'e}partement de Math{\'e}matiques \\
Universit{\'e} du Qu{\'e}bec
{\`a} Montr{\'e}al \\
Case postale 8888 \\ succursale centre--ville\\ Montr{\'e}al (Qu{\'e}bec)\\
H3C 3P8, Canada} \email{apostolo@math.uqam.ca}

\address{Tedi Dr\u{a}ghici \\ Department of Mathematics \\ Florida
International University \\ Miami FL 33199 \\ USA}
\email{draghici@fiu.edu}

\address{Andrei Moroianu \\ CMAT\\ {\'E}cole Polytechnique \\
UMR 7640 du CNRS\\ 91128 Palaiseau \\ France}
\email{am@math.polytechnique.fr}

\begin{document}
\begin{abstract}

An odd--dimensional version of the Goldberg conjecture
was formulated and proved in \cite{bg}, by using an orbifold analogue of
Sekigawa's arguments in \cite{se}, and an approximation argument  of
K--contact structures with quasi--regular ones. We provide here another
proof of this result.

\vspace{0.1cm}
\end{abstract}

\maketitle

\section{Introduction}

The celebrated Goldberg conjecture states that every compact almost
K{\"a}hler Einstein manifold $M$ is actually K{\"a}hler--Einstein. This
conjecture  was confirmed by Sekigawa \cite{se} in the case when $M$ has
non--negative scalar curvature. The odd--dimensional analogues of
K{\"a}hler manifolds are Sasakian manifolds, and those of almost
K{\"a}hler manifolds are
K--contact manifolds. In \cite{bg}, Boyer and Galicki proved the following
odd--dimensional analogue of Goldberg's conjecture.

\begin{ath}\label{main}{\em \cite{bg}} Any compact Einstein K--contact
manifold $(M,g,\xi)$ is Sasakian.
\end{ath}

Their proof goes roughly as follows. First, an Einstein K--contact
manifold has prescribed (positive) Einstein constant. If the
K--contact structure is quasi--regular (i.e. the orbits of the Reeb vector
field $\xi$ are closed), then the quotient of $M$ by the flow of $\xi$ is
an almost K{\"a}hler orbifold \cite{thomas} which is
Einstein with positive scalar curvature by the O'Neill formulas. One then
applies Sekigawa's proof to obtain that the almost K{\"a}hler structure is
integrable, which in turn means that the  K--contact structure is
Sasakian.  If the K--contact structure is not
quasi--regular, the space of orbits of $\xi$ is not an orbifold
(and may not be even a tractable topological space). To overcome this
difficulty,  the authors of
\cite{bg} provide a beautiful argument showing that the Reeb vector field
$\xi$ can be approximated (in a suitable sense) by a sequence of
quasi--regular Reeb vector fields
$\xi_i$ which define K--contact structures on a sequence of (no longer
Einstein) metrics
$g_i$ approaching $g$. Then for the sequence of orbifolds thus
obtained, one can use ``approximative'' Sekigawa formulas and
eventually show that the K--contact structure is integrable.

The aim of this note is to give another proof of Theorem
\ref{main} and to study further possible extensions. 
Instead of the quotient of $M$ by the Reeb flow, we
consider another almost K{\"a}hler manifold naturally associated
to $M$, namely the cone over $M$. It is well--known that the cone
is a smooth, {\it non--compact} Ricci--flat almost K{\"a}hler
manifold which is K{\"a}hler if and only if $M$ is Sasakian. It
therefore suffices to prove the integrability of the almost
K{\"a}hler cone structure. It would seem to be difficult to apply
directly Sekigawa's arguments in this situation because of the
non--compactness of the cone. But this can be overcome easily: we
first apply a point--wise version of Sekigawa's formula on the cone
manifold, and then integrate it on the level sets of the radial
function (which are compact manifolds).

The use of this approach tempted us to extend the conjecture to
the more general case of {\it contact metric structures}, when the
metric is no longer bundle--like. Indeed, one could argue that the
analogue of almost K{\"a}hler manifolds in odd dimensions are the
contact metric structures, since they correspond to the level sets
of the radial function of almost K{\"a}hler cone metrics. The
contact analogue of the Goldberg conjecture would then assert that
{\it any compact Einstein contact metric manifold is
Sasakian--Einstein}. This statement turns out to be false in
general, as it follows from an example of D. Blair on the flat
3--torus, which we recall in the last section. Still, the
counterexample does not generalize to higher dimensions, so the
problem seems to be worth further investigation. We make a 
step in this direction; in the particular case when the Einstein
metric admits already a compatible Sasakian structure, we use
Theorem~\ref{main} to show that any other compatible contact
metric structure is necessarily Sasakian.
\begin{ath}\label{corollary}
Let $(M,g,\xi)$ be a compact Sasakian--Einstein manifold of
dimension $2n+1$. Then any contact metric structure $(\xi',g)$ on
$(M,g)$ is Sasakian.  Moreover, if  $\xi'$ is different from $\pm
\xi$, then the following two cases occur:
\begin{enumerate}
\item[$\bullet$] $(M,g)$ admits a 3--Sasakian structure and $\xi$
and $\xi'$ belong to the underlying $S^2$--family of  Sasakian structures.
\item[$\bullet$] $(M,g)$ is covered by the round sphere $S^{2n+1}$.
\end{enumerate}
\end{ath}
Note that the cone construction identifies the set  of all
Sasakian structures on the round sphere $S^{2n+1}$ with the
homogeneous space $SO(2n+2)/U(n+1)$.

\section{Preliminaries}

Let $(M,g)$ be a Riemannian manifold. We define the cone
$\bar M:=M\x\RM^*_+$ endowed with the metric $\gg=dr^2+r^2g$, and denote
by
$\nn$ the covariant derivative of $\gg$. It is well--known that the cone is
a non--complete Riemannian manifold which can be completed at $r=0$ if and
only if $M$ is a round sphere.

Every vector field $X$ on $M$ induces in a canonical way a vector
field $(X,0)$ on $\bar M$, which (with a slight abuse of notation) will
still be denoted by $X$. Similarly, we denote by the same symbol the forms
on $M$ and their pull--backs to $\bar M$ (with respect to the projection on
the first factor). Let us denote by $\dr$ the vector field
$\frac{\pa}{\pa r}$ on $\bar M$. The following formulas relate the
covariant derivatives $\n$ and $\nn$, and are immediate consequences of
the definitions.
\begin{equation}\label{vec}
\nn_\dr\dr=0; \quad \nn_X\dr=\nn_\dr X=\frac{1}{r}X;\quad
\nn_XY=\n_XY-rg(X,Y)\dr.
\end{equation}
Using this, we obtain for every vector $X$ and a $p$--form $\o$ on $M$
\begin{equation}\label{for}
\nn_\dr\o=-\frac{p}{r}\o\quad\hbox{and}\quad\nn_X\o=\n_X\o-
\frac{1}{r}dr\wedge X\i\o,
\end{equation}

\begin{equation}\label{for2}
\nn_\dr dr=0\quad\hbox{and}\quad\nn_Xdr=rX^\flat.
\end{equation}

The curvature tensors $R$ and $\bar R$ of $M$ and $\bar M$, respectively,
are related by
\begin{equation}\label{cou}
\bar R(\dr,\cd)=0\quad\hbox{and}\quad \bar
R(X,Y)Z=R(X,Y)Z+g(X,Z)Y-g(Y,Z)X.
\end{equation}
\begin{ede}  A {\em contact metric structure} on a Riemannian
manifold
$M$ is a unit length vector field
    $\xi$ such that the 1--form $\eta:=\la\xi,\cdot\ra$ and the
    endomorphism $\f$ associated to
$\frac{1}{2}d\eta$ are inter--related by
\begin{equation}
\label{kc}\f^2=-1+\eta\otimes\xi.
\end{equation}
Since $\f^2(\xi)=0$, we get $|\f (\xi)|^2=-\la\xi,\f^2(\xi)\ra=0$,
so $\f(\xi)=0$. In other words, $\f$ defines a complex structure on
the distribution orthogonal to $\xi$.

A contact metric structure $(M,g,\xi,\f,\eta)$ is called {\em
K--contact} if $\xi$ is Killing.
The  K--contact structure $(M,g,\xi,\f,\eta)$ is called {\em
Sasakian} if
\begin{equation}\label{xd}{\nabla_\cdot}\n\xi=\xi\wedge
\cdot.\end{equation}
\end{ede}

Given a contact metric manifold $(M,g,\xi,\f,\eta)$, we construct
a 2--form $\O$ on $\bar M$, defined by
\begin{equation}\label{om}
\O=r dr\wedge\eta+\frac{r^2}{2}d\eta.
\end{equation}
This 2--form is clearly compatible with $\gg$, and therefore defines an
almost complex structure $J$ on $\bar M$ by
$\Omega(\cd,\cd)=\gg(J\cd,\cd)$. Moreover, $\O$ is obviously closed,
meaning that $(\bar M,J)$ is almost K{\"a}hler.
It is well--known that $\Omega$ is parallel (i.e. $(\bar M,J)$ is
K{\"a}hler) if and only if the contact structure $\xi$ is Sasakian.

We close this section with the following

\begin{elem}\label{le1} {\rm (i)} The codifferentials on $M$ and $\bar M$
are related by
\begin{equation}
\d^{\bar M}(r^k\s)=r^{k-2}\d^M\s,\quad \forall\s\in\L^1M.
\end{equation}

\noindent
{\rm (ii)} The Laplacians on $M$ and $\bar M$ are related by
\begin{equation}
\Delta^{\bar M}(r^kf)=r^{k-2}(\Delta^Mf-k(2n+k)f), \quad\forall f\in
C^\infty(M).
\end{equation}
\end{elem}

\begin{proof} (i) If $(e_i)$ denotes a local
orthonormal base on $M$, we have
\beq
\d^{\bar
M}(r^k\s)&=&\sum_i(-\ei(r^k\s(\ei))+r^k\s(\nn_{\ei}\ei))-\dr(r^k\s(\dr))
+r^k\s(\nn_{\dr}\dr)\\
&=&\sum_i-r^{k-2}e_i(\s(e_i))+r^{k-2}\s(\n_{e_i}e_i)-r\dr)=r^{k-2}\d^M\s.
\eeq

(ii) Similarly,
\beq\Delta^{\bar
M}(r^kf)&=&\sum_i(-\ei(\ei(r^kf))+\nn_{\ei}\ei(r^kf))-\dr(\dr(r^kf))\\
&=&\sum_i(-r^{k-2}e_i(e_i(f))+\frac{1}{r^2}
(\n_{e_i}e_i-r\dr)(r^kf))-k(k-1)r^{k-2}f\\
&=&r^{k-2}\Delta^Mf-k(2n+1)r^{k-2}f-k(k-1)r^{k-2}f\\
&=&r^{k-2}(\Delta^Mf-k(2n+k)f). \;
\eeq
\r

\section{Proof of Theorem \ref{main}}
Let $(M^{2n+1},g,\xi)$ be a compact K--contact Einstein manifold.
By a result of Blair (\cite{bl}, Theorem 7.1), a contact metric
manifold is K--contact if and only if $\Ric(\xi,\xi) = 2n$; thus,
the Einstein constant in our case must be $2n$.

Consider now the cone $\bar M$, which is an almost K{\"a}hler
manifold. We use the following Weitzenb{\"o}ck--type formula, taken from
\cite[Prop.2.1]{adm}.

\begin{epr} \label{prop2} For any almost K{\"a}hler manifold
$(\bar M,\gg, J, \O)$ with covariant derivative denoted by $\nn$ and
curvature tensor $\bar R$, the following point--wise
relation holds:

\begin{eqnarray} \label{la}
\Delta(s^* - s) &=& -4 \delta( J \delta^{\nn} (J\bar\Ric'')) +
8 \delta (\langle \bar\rho^* , \nn_{\cdot} \; \O\rangle) +
2|\bar Ric''|^2\\ \nonumber
& &  - 8|\bar R''|^2  - |\nn^*\nn \O |^2 - |\phi|^2
+ 4\langle {\rho}, \phi \rangle - 4\langle {\rho}, \nn^*\nn \O \rangle
\; ,
\end{eqnarray}
where: $\phi(X,Y) = \langle \nn_{JX} \O,\nn_Y \O \rangle $, $s$
and $s^*$ are respectively the scalar and $*$--scalar curvature,
$\bar\Ric''$ is the $J$--anti--invariant part of the Ricci tensor
$\bar\Ric$, $\rho$ is the $(1,1)$--form associated to the
$J$--invariant part of $\bar\Ric$, $\bar\rho^*:=\bar R(\O)$ and
$\bar R''$ denotes a certain component of the curvature tensor.
\end{epr}

In our situation, since $M^{2n+1}$ is Einstein with constant $2n$,
(\ref{cou}) shows that $\bar M$ is Ricci--flat. So the formula above
becomes
\begin{equation}\label{la1}
\Delta^{\bar M} s^*-
8 \delta^{\bar M} (\langle \bar\rho^* , \nn_{\cdot} \; \O\rangle) = -
8|\bar R''|^2  - |\nn^*\nn \O |^2 - |\phi|^2
\end{equation}
We now use Lemma \ref{le1} in order to express the left--hand side of this
equality in terms of the codifferential and Laplacian on $M$.
  From (\ref{cou}) we get $\bar\rho^*(X,\dr)=0$ and
$\bar\rho^*(X,Y)=\bar g(\bar R(\ei,J\ei)X,Y)= \rho^*(X,Y)$, for some
2--form $\rho^*$ on $M$. Taking the scalar product with $\O$ yields
\begin{equation}\label{op}
s^*=\frac{1}{r^2}f
\end{equation}
for some function $f$ on $M$. Note that $f$ is everywhere
positive on $M$ since $s^*=s^*-s=|\nn\O|^2$ on $\bar M$ (see e.g.
\cite{adm}, p.~777).

Now, from (\ref{for}), (\ref{for2}) and (\ref{om}) we get
$\nn_\dr\O=0$ and $\nn_X\O=r^2\o+rdr\wedge\tau_X$ for some 2--form $\o$
and 1--form $\tau_X$ on $M$. Consequently, the 1--form $\langle
\bar\rho^* , \nn_{\cdot} \; \O\rangle$ on $\bar M$ is easily seen to
be of the form
\begin{equation}\label{op2}
\langle\bar\rho^* , \nn_{\cdot} \; \O\rangle=\frac{1}{r^2}\a
\end{equation}
for some 1--form $\a$ on $M$. Using (\ref{op}), (\ref{op2}) and
Lemma~\ref{le1}, the equality (\ref{la1}) becomes
\begin{equation}\label{la2}
\frac{1}{r^4}(\Delta^Mf+2(2n-2)f-8\d^M\a)= -
8|\bar R''|^2  - |\nn^*\nn \O |^2 - |\phi|^2.
\end{equation}
Integrating this last equation on each level set $M_r:=\{r=constant\}$ of
$\bar M$ yields
$$\int_{M_r}\frac{2(2n-2)}{r^4}f+8|\bar R''|^2 +|\nn^*\nn \O |^2
+|\phi|^2.$$
In particular, since $f\ge 0$, $\phi$ vanishes identically on $\bar
M$, hence $|\nn_X\O|^2=-\phi(X,JX)=0$ for every $X$ on $\bar M$. Thus
$\bar M$ is K{\"a}hler, so $M$ is Sasakian.

\section{Proof of Theorem \ref{corollary}}
Let $(M,g,\xi)$ be a compact Sasakian--Einstein manifold and
$\xi'$ be another $g$--compatible contact metric structure. Since
${\rm Ric} = 2n g$, it follows that $(g, \xi')$ is ${\rm
K}$--contact (see \cite{bl}, Theorem 7.1), hence Sasakian
according to Theorem~\ref{main}. Since $(M,g)$ is complete, by a
result of Gallot \cite{gallot}, the cone $({\bar M}, {\bar g})$
over $(M,g)$ is (locally) de~Rham irreducible unless it is flat
(i.e. $(M,g)$ is of positive constant curvature). Therefore,
Theorem~\ref{corollary}  follows from the  following general
observation.
\begin{epr}\label{hypersasaki}
Suppose $(M,g)$ is a Riemannian manifold whose cone is locally
irreducible and which admits two Sasakian structures  $(g,\xi)$
and $(g,\xi')$ with $\xi \neq \pm \xi'$.  Then $(M,g)$ admits a
3--Sasakian structure and $\xi$ and $\xi'$ belong to the
underlying $S^2$--family of Sasakian structures.
\end{epr}
\begin{proof}
Let $({\bar M}, {\bar g})$ be the cone over $(M,g)$. The Sasakian
structures $\xi$ and $\xi'$ give rise to two K{\"a}hler structures,
$J$ and $J'$,  on $(\bar M, \bar g)$ with $J \neq \pm J'$ (because
$\xi \neq \xi'$ by assumption).  It suffices to show that $(\bar
M, \bar g)$ must be hyperk{\"a}hler and $J$ and $J'$ belong to the
$S^2$--family of compatible K{\"a}hler structures. The
anti--commutator $Q= JJ' + J'J$ of $J$ and $J'$ is symmetric and
parallel with respect to ${\bar g}$; since  $({\bar M}, \bar g )$
is locally irreducible, $Q= \lambda {\rm Id}$ for some real
constant $\lambda$. Since $J$ and $J'$ are both $\bar g$
orthogonal, the Cauchy--Schwartz inequality implies $|\lambda | \le
2$; it is easy to see that equality is possible if only if $J= \pm
J'$, a situation that we excluded. Similarly, the commutator $A=
JJ' - JJ'$ of $J$ and $J'$ is parallel and skew--symmetric with
respect to the metric ${\bar g}$ and by using the corresponding
property of $Q$, it verifies $A^2 = (\lambda^2 - 4) {\rm Id}$. It
follows that $I = \frac{1}{\sqrt {4 - \lambda^2}} A$ defines a
parallel, ${\bar g}$--compatible complex structure on $(\bar M,
\bar g)$ which anti--commutes with both $J$ and $J'$; therefore
$(\bar g, I, J, K= IJ)$ defines a hyperk{\"a}hler structure. The
equality  $JJ' + J'J = \lambda {\rm Id}$ also shows that $J'$
belongs to the $S^2$--family of K{\"a}hler structures generated by
$(\bar g, I, J, K)$. \end{proof}

\section{An example and further comments}

As explained in the introduction, it was tempting to ask the
following question, slightly more general than Theorem \ref{main}:
{\it is every compact Einstein contact metric manifold
Sasakian--Einstein?} The answer is negative in general, as the
following simple example of D. Blair shows (see \cite{bl}, p. 23,
p. 68--69 \& p. 52--53).
\begin{exa} The 1--form
$\eta:=\cos t\, dx+\sin t\, dy$ defines a (non--regular) contact
metric structure on the flat torus $T^3$ {\rm (}where $t,\ x$ and
$y$ are standard coordinates on $T^3$ of periods  $2\pi${\rm )},
which is not K--contact (and hence not Sasakian).
\end{exa}

Note however that this is the only negative example to the above
question in dimension 3. Indeed, in this dimension, Blair and
Sharma \cite{blair-sharma} proved that a contact metric manifold
of constant curvature has either curvature $+1$ and is Sasakian,
or curvature $0$ and is isometric to the above example. Note also
that the example does not directly generalize to higher
dimensions, as Blair \cite{blair-tohoku} also shows that there are
no flat contact metric structures in dimension $ \geq 5$.  More
generally, there is a theorem of Olszak \cite{olszak} that in
dimension $ \geq 5$ there are no contact metric manifold of
constant curvature, unless the curvature is $+1$ and the structure
Sasakian. Hence, there are reasons to still investigate the above
question.

\vspace{0.1cm} \noindent {\bf Acknowledgments:} The authors are
grateful to David Blair, Charles Boyer and Lieven Vanhecke  for
useful comments.

\labelsep .5cm

\end{document}